\documentclass[10pt]{article}

	\usepackage[latin1]{inputenc}
	\usepackage{graphicx, amsmath, amsthm, amsfonts, amssymb, subfigure, amstext}
	\usepackage{wrapfig, amsopn, latexsym, epsfig, multirow, verbatim}
	\usepackage{comment, caption, bm}
	\usepackage[plainpages=false,pdfpagelabels,hypertexnames,linktocpage,
		colorlinks, citecolor=red,linkcolor=blue]{hyperref}
	\usepackage[numbers]{natbib}
	\usepackage{microtype}

	\newtheorem{theorem}{Theorem}[section]
	\newtheorem{remark}{Remark}[section]

	\usepackage[hmargin=2cm,vmargin={3.5cm,4cm}]{geometry}
	\numberwithin{equation}{section}
	\allowdisplaybreaks

	\title{Transient behavior of fractional queues and related processes}

	\author{$\text{Dexter O. Cahoy}_1$, $\text{Federico Polito}_2^{*}$, $\text{Vir Phoha}_3$\\
		\footnotesize (1) -- Department of Mathematics and Statistics\\
	    \footnotesize College of Engineering and Science, Louisiana Tech University, USA\\
	    \footnotesize Tel: +1 318 257 3529, fax: +1 318 257 2182\\
		\footnotesize Email address: dcahoy@latech.edu\\
		\footnotesize (2) -- Department of Mathematics, University of Torino, Italy\\
		\footnotesize Tel: +39 011 670 2937, fax: +39 011 670 2878\\
		\footnotesize Email address: federico.polito@unito.it\\
		\footnotesize (3) -- Department of Computer Science\\
        \footnotesize College of Engineering and Science, Louisiana Tech University, USA\\
		\footnotesize Tel: +1 318 257 2298, fax: +1 318 257 4922\\
		\footnotesize Email address: phoha@latech.edu\\
		}
	
\begin{document}

		\maketitle

		\begin{abstract}
		
			We propose a generalization of the classical M/M/1 queue process.
			The resulting model is derived  by applying fractional derivative operators to a system of
			difference-differential equations.  This generalization includes
			both non-Markovian and Markovian properties which naturally provide greater flexibility
			in modeling real queue systems than its  classical counterpart.  Algorithms to simulate M/M/1
			queue process  and the related  linear birth-death process are provided.  Closed-form expressions
			of the point and interval estimators of  the parameters of the proposed fractional stochastic models are
			also presented. These methods are necessary to make these models usable in practice.
			The proposed fractional M/M/1 queue model and the statistical methods are illustrated using financial data.
			 
			\vspace{0.1in}

		\end{abstract}

		\section{Introduction}

			The M/M/1 queue is without a doubt the simplest model for a queue process.
			It is characterized by arrivals determined by a Poisson process and an independent service time which
			is negative-exponentially distributed.  It is relatively simple and yet 
			the analysis of its transient behavior leads to considerable difficulties.
			The main source of these difficulties is the presence of a non-absorbing boundary at zero
			(empty queue).  This means that the analysis becomes simpler when we consider models with absorbing boundaries.  As a direct result,   the state probability of 
			a linear birth-death process, that is the probability that
			the queue length is  $n$ at a specific time $t$, has a particularly
			nice form. 

			The aim of this paper is to study some related point processes governed
			by difference-differential equations containing fractional derivative operators.
			These processes	are direct generalizations of the classical M/M/1 queue
			and the linear birth-death processes.
			It is well-known that a fractional derivative operator induces a non-Markovian behavior
			into a system \citep[see][]{vetaq}.  Moreover,   parameter estimation  and path generation algorithms of these new fractional stochastic models are derived.   Note that the proposed fractional point models (with Markovian and non-Markovian properties) are parsimonious which makes them desirable for modeling real-world non-Markovian queueing systems.   Observe that fractional point processes driven by fractional difference-differential equations  such as the fractional Poisson,  the  fractional birth, the fractional death,  and the fractional birth-death processes have already been gaining attention more  recently \citep[see, e.g.,][]{laskin, beg, cuw10, polito, fedo}. 
			
			The article is structured as follows.
			Section \ref{secsec} presents the explicit construction of the fractional M/M/1 queue
			starting from the governing equations and a particular subordination relation.
			The main result derived in this section is the explicit form of the
			transient state probabilities for each value of the parameter of fractionality.
			Information regarding the steady-state behavior (stationary behavior) is also highlighted; in
			particular the fractional process shares the same steady-state behavior
			as the classical non-fractional case.
			In Section \ref{secsec2} we develop closed-form estimators (point and interval)
			for the model parameters in the case of the fractional linear birth-death process. This is preparatory
			for the similar subsequent analysis applied to the M/M/1 queue (Section \ref{secsec3}).
			The article ends with an application which shows  that our constructed estimators perform
			well in a real-world example.

		\section{Results for a fractional process related to M/M/1 queues}

			\label{secsec}
			The classical M/M/1 queue process $N(t)$, $t \ge 0$, that is the queue length
			in time can be described by the following difference-differential
			equations governing the state probabilities $p_k(t) = \Pr \{ N(t) = k | N(0) = i \}$, $k \ge 0$: 
			\begin{align}
				\label{20xxx}
				\begin{cases}
					\frac{\mathrm d}{\mathrm dt} p_k(t) = -(\lambda+\mu) p_k(t) + \lambda p_{k-1}(t)
					+ \mu p_{k+1}(t), & k \ge 1, \\
					\frac{\mathrm d}{\mathrm dt} p_0(t) = - \lambda p_0(t) +\mu p_1(t), \\
					p_k(0) = \delta_{k,i},
				\end{cases}
			\end{align}
			where $i \in \mathbb{N} \cup \{0\}$ is the initial number of individuals in the queue and
			$\delta_{k,i}$ is the Kronecker's delta. In \eqref{20xxx} $\lambda>0$ and $\mu>0$ are 
			the entrance and the service rates, respectively.
			
			To arrive at a possible fractional model
			we consider the Caputo fractional derivative $D^\alpha_t$, $\alpha \in (0,1]$, with respect to time $t$.
			If $p_k^\alpha(t) = \Pr \{ N^\alpha(t) = k \}$, $k \ge 0$,
			where $N^\alpha(t)$,
			is the fractional $M/M/1$ queue with parameter $\alpha$,
			the generalized difference-differential equations for the state probabilities
			with arrival rate $\lambda>0$, service rate $\mu>0$ and $i \ge 0$ initial customers, read
			\begin{align}
				\label{20xx}
				\begin{cases}
					D^{\alpha}_t p_k^\alpha(t) = -(\lambda+\mu) p_k^\alpha(t) + \lambda p_{k-1}^\alpha(t)
					+ \mu p_{k+1}^\alpha(t), & k \ge 1, \\
					D^{\alpha}_t p_0^\alpha(t) = - \lambda p_0^\alpha(t) +\mu p_1^\alpha(t), \\
					p_k^\alpha(0) = \delta_{k,i}.
				\end{cases}
			\end{align}
			First, we will follow \citet{bailey1954continuous,bailey1990elements}
			for the derivation of the probabilities $p_k^\alpha(t)$, $k \ge 0$, $t \ge 0$
			but adapting the method to take into considerations the presence of the Caputo derivative.
			The result obtained by Bailey is the so-called classical solution in terms
			of modified Bessel functions of the first kind.  
			Note however that the derivation of the state probabilities in the classical case $\alpha=1$
			can be carried out in several equivalent ways (see for example \citet{champer, partha, abate}).
			In the following we will first treat the solution derived by Bailey and then we will use
			a simpler but lesser known form due to \citet{sharma}.
			
			We
			indicate  $G^\alpha(z,t) = \sum_{k=0}^\infty z^k p_k^\alpha(t)$ as the probability
			generating function.
			
			\begin{theorem}			
			
				The Laplace transform $\tilde{G}^\alpha(z,s) = \int_0^\infty e^{-st}
				G^\alpha(z,t) \, \mathrm dt$, $\alpha \in (0,1]$, can be written as
				\begin{align}
					\label{gener}
					\tilde{G}^\alpha(z,s) = s^{\alpha-1} \frac{z^{i+1}-(1-z)\left[ a_2(s) \right]^{i+1}
					\left[ 1-a_2(s) \right]^{-1}}{-\lambda \left[ z-a_1(s) \right] \left[ z-a_2(s) \right]},
					\qquad |z| \le 1, \: \Re(s) > 0.
				\end{align}
				where $a_1(s)$ and $a_2(s)$ are the zeros of $f(z,s) = zs^\alpha-(1-z)(\mu-\lambda z)$.
	
				\begin{proof}		

					From \eqref{20xx}, we can write
					\begin{align}
						D^\alpha_t \left[ G^\alpha(z,t)-p_0^\alpha(t) \right] = -(\lambda+\mu)
						\left[ G^\alpha(z,t)-p_0^\alpha(t) \right]
						+\lambda z G^\alpha(z,t).
					\end{align}
					Using the equation on $p_0^\alpha(t)$ we have
					\begin{align}
						D^\alpha_t G^\alpha(z,t) & = -\lambda G^\alpha(z,t) -\mu G^\alpha(z,t) +\mu p_0^\alpha(t)
						+\lambda z G^\alpha(z,t) +\frac{\mu}{z}
						\left[ G^\alpha(z,t)-p_0^\alpha(t) \right],
					\end{align}
					and after simplifying, we obtain, for $|z|\le 1$, the Cauchy problem
					\begin{align}
						\label{burburu}
						\begin{cases}
							z D^\alpha_t G^\alpha(z,t) = (1-z) \left[ G^\alpha(z,t) (\mu-\lambda z)
							-\mu p_0^\alpha(t) \right], \\
							G^\alpha(z,0) = z^i.
						\end{cases}
					\end{align}
					Applying the Laplace transform $\tilde{G}^\alpha(z,s)
					= \int_0^\infty e^{-st} G^\alpha(z,t) \, \mathrm dt$ to \eqref{burburu} leads to
					\begin{align}
						z \left[ s^\alpha \tilde{G}^\alpha(z,s) - s^{\alpha-1} G^\alpha(z,0) \right]
						= (1-z) \left[ \tilde{G}^\alpha(z,s)
						(\mu-\lambda z) - \mu \tilde{p}_0^\alpha(s) \right],
					\end{align}
					where $\tilde{p}_0^\alpha(s) = \int_0^\infty e^{-st} p_0^\alpha(t) \, \mathrm dt$.
					After some simple algebraic calculations
					we then have
					\begin{align}
						\label{00aa}
						\tilde{G}^\alpha(z,s) = \frac{s^{\alpha-1} z^{i+1}-\mu(1-z)
						\tilde{p}_0^\alpha(s)}{zs^\alpha-(1-z)(\mu-\lambda z)},
						\qquad |z| \le 1, \: \Re(s) > 0.
					\end{align}
					As the above function converges in $|z|\le 1$, the zeros of the numerator and the denominator should
					coincide. Let us indicate the zeros of the numerator as
					\begin{align}
						a_{12} (s) = \frac{s^\alpha +\lambda+\mu \pm \left[ (s^\alpha +\lambda +\mu)^2 -4\lambda
						\mu \right]^{1/2}}{2\lambda},
					\end{align}
					with $|a_2(s)|<|a_1(s)|$, $\Re(s) > 0$. Note that
					\begin{align}
						\label{rel}
						\begin{cases}
							a_1(s)+a_2(s) = (s^\alpha + \lambda +\mu)/\lambda, \\
							a_1(s) a_2(s) = \mu/\lambda, \\
							-\lambda[1-a_2(s)][1-a_1(s)] = s^\alpha.
						\end{cases}
					\end{align}
					By Rouch\'e theorem \citep[Page 168]{krantz} we have that the only zero in the unit circle
					is $a_2(s)$. Therefore it follows that
					\begin{align}
						s^{\alpha-1} \left[ a_2(s) \right]^{i+1} - \mu [1-a_2(s)] \tilde{p}_0^\alpha(s) = 0,
					\end{align}
					which gives
					\begin{align}
						\label{pirelli}
						\tilde{p}_0^\alpha(s) = \frac{s^{\alpha-1}\left[ a_2(s) \right]^{i+1}}{\mu[1-a_2(s)]}.
					\end{align}
					Now, by considering that
					\begin{align}
						zs^\alpha-(1-z)(\mu-\lambda z) = -\lambda \left[ z-a_1(s) \right] \left[ z-a_2(s) \right],
					\end{align}
					equation \eqref{00aa} can be rewritten as
					\begin{align}
						\label{gener2}
						\tilde{G}^\alpha(z,s) = s^{\alpha-1} \frac{z^{i+1}-(1-z)\left[ a_2(s) \right]^{i+1}
						\left[ 1-a_2(s) \right]^{-1}}{-\lambda \left[ z-a_1(s) \right] \left[ z-a_2(s) \right]}.
					\end{align}
			
				\end{proof}

			\end{theorem}			
			
			In the following Theorem \ref{teo1} we prove a
			subordination relation for the fractional queue $N^\alpha(t)$, $t \ge 0$, $\alpha \in (0,1]$.
			This is essential for our next results.
			Before that, let us introduce some facts on the $\alpha$-stable subordinator and its inverse process.
		
			Let us call $V^\alpha(t)$, $t \ge 0$, the $\alpha$-stable subordinator (see for details \citet{bertoin},
			cap. III) and
			let us define its inverse process as its hitting time
			\begin{align}
				\label{54}
				E^\alpha(t) = \inf \{ s > 0 \colon V^\alpha(s) > t \}. 
			\end{align}
			The processes $V^\alpha(t)$ and $E^\alpha(t)$ are characterized by their Laplace transforms.
			For the $\alpha$-stable subordinator we have
			\begin{align}
				\mathbb{E}e^{-\xi V^\alpha_t} = e^{-t \xi^\alpha}, \qquad \alpha \in (0,1], % \: t \ge 0,
			\end{align}
			and for its inverse process the time-Laplace transform reads
			\begin{align}
				\int_0^\infty e^{-\xi t} \left( \Pr \{ E^\alpha(t) \in \mathrm ds \}/\mathrm ds \right)
				\mathrm dt = \xi^{\alpha-1} e^{-s \xi^\alpha} \qquad \alpha \in (0,1]. %\: t \ge 0.
			\end{align}
		
			\begin{theorem}
				\label{teo1}
				Let $N^1(t)=N(t)$, $t \ge 0$, be the classical $M/M/1$ queue and let
				$E^\alpha(t)$, $t \ge 0$, $\alpha \in (0,1]$, be an inverse $\alpha$-stable subordinator \eqref{54}
				independent of $N^1(t)$.
				The fractional $M/M/1$ queue $N^\alpha(t)$, $t \ge 0$, $\alpha \in (0,1]$, can be represented as
				\begin{align}
					\label{20sub}
					N^\alpha(t) = N^1(E^\alpha(t)), \qquad t \ge 0, \: \alpha \in (0,1),
				\end{align}
				where the equality holds for the one-dimensional distribution.
			
				\begin{proof}
					Let us consider the initial value problem
					\begin{align}
						\begin{cases}
							z D^\alpha_t G^\alpha(z,t) = (1-z) \left[ G^\alpha(z,t) (\mu-\lambda z)
							-\mu p_0^\alpha(t) \right], \\
							G^\alpha(z,0) = z^i,
						\end{cases}
					\end{align}
					which is equivalent to \eqref{20xx}. Applying the Laplace transform we obtain
					\begin{align}
						\label{20ff}
						z \left[ s^\alpha \tilde{G}^\alpha(z,s) - s^{\alpha-1} G^\alpha(z,0) \right]
						= (1-z) \left[ \tilde{G}^\alpha(z,s)
						(\mu-\lambda z) - \mu \tilde{p}_0^\alpha(s) \right],
					\end{align}
					Note that if \eqref{20sub} holds we can write
					\begin{align}
						\label{20dd}
						\tilde{G}^\alpha(z,s) & = \int_0^\infty e^{-st} \left[ \sum_{k=0}^\infty
						z^k \int_0^\infty \Pr \{ N^\alpha(y) = k \} \Pr \{ E^\alpha(t) \in \mathrm dy \} \right]
						\mathrm dt \\
						& = \int_0^\infty e^{-st} \left[ \int_0^\infty G(z,y) \Pr \{ E^\alpha(t) \in \mathrm dy \}
						\right] \mathrm dt \notag \\
						& = \int_0^\infty G(z,y) s^{\alpha-1} e^{-ys^\alpha} \mathrm dy, \notag
					\end{align}
					and 
					\begin{align}
						\label{20ee}
						\tilde{p}_0^\alpha(s) & = \int_0^\infty e^{-st} p_0^\alpha(t) \mathrm dt \\
						& = \int_0^\infty e^{-st} \left[ \int_0^\infty p_0(y) \Pr \{ E^\alpha(t) \in \mathrm dy \} \right]
						\mathrm dt \notag \\
						& = \int_0^\infty p_0(y) s^{\alpha-1} e^{-ys^\alpha} \mathrm dy. \notag
					\end{align}
					We now show that \eqref{20dd} and \eqref{20ee} satisfy \eqref{20ff}. Observe that
					\begin{align}
						\label{20gg}
						& z \left[ s^\alpha \int_0^\infty G(z,y) e^{-ys^\alpha} \mathrm dy \mu
						- z^i \right] 
						 = (1-z) \left[ (\mu-\lambda z) \int_0^\infty G(z,y) e^{-ys^\alpha} \mathrm dy
						- \mu \int_0^\infty p_0(y) e^{-ys^\alpha} \mathrm dy \right].
					\end{align}
					Consider the right hand side of \eqref{20dd}. We can write
					\begin{align}
						& (1-z) \left[ (\mu-\lambda z) \int_0^\infty G(z,y) e^{-ys^\alpha} \mathrm dy
						- \mu \int_0^\infty p_0(y) e^{-ys^\alpha} \mathrm dy \right] \\
						& = \int_0^\infty e^{-ys^\alpha} (1-z) \left[ (\mu-\lambda z) G(z,y)
						- \mu p_0(y) \right] \mathrm dy. \notag
					\end{align}
					Considering that $G(z,y)$ and $p_0(y)$ satisfy
					\begin{align}
						z \frac{\partial}{\partial y} G(z,y) = (1-z) \left[ (\mu-\lambda z) G(z,y)
						-\mu p_0(y) \right],
					\end{align}
					we immediately obtain that
					\begin{align}
						& (1-z) \left[ (\mu-\lambda z) \int_0^\infty G(z,y) e^{-ys^\alpha} \mathrm dy
						- \mu \int_0^\infty p_0(y) e^{-ys^\alpha} \mathrm dy \right] \\
						& = z \int_0^\infty e^{-y s^\alpha} \frac{\partial}{\partial y} G(z,y) \mathrm dy \notag \\
						& = z \left[ \left. G(z,y) e^{-ys^\alpha} \right|_{y=0}^{y=\infty}
						+ s^\alpha \int_0^\infty G(z,y) e^{-y s^\alpha} \mathrm dy \right] \notag \\
						& = z \left[ s^\alpha \int_0^\infty G(z,y) e^{-ys^\alpha} \mathrm dy
						- z^i \right]. \notag
					\end{align}
					This concludes the proof.
				\end{proof}
			\end{theorem}
			
			Using the Laplace transform \eqref{gener} and the calculations
			carried out in \citet{bailey1990elements}
			we can gain some insights on the mean value of the process.
			
			\begin{theorem}			
				
				We have that
				\begin{align}
					\label{meap}
					\mathbb{E} N^\alpha(t) =  i + (\lambda-\mu) \frac{t^\alpha}{\Gamma(\alpha+1)}
					+ \mu J^\alpha p_0^\alpha(t),
				\end{align}
				where
				\begin{align}
					J^\alpha f(t) = \frac{1}{\Gamma(\alpha)} \int_0^t (y-t)^{\alpha-1} f(y) \mathrm dy, \qquad t > 0,
				\end{align}
				is the Riemann--Liouville fractional integral \citep{samko}.
			
				\begin{proof}
			
					By means of the Laplace transform \eqref{gener} of the probability generating
					function we can write
					\begin{align}
						\label{lap-mea}
						\tilde{\mathbb{E}} N^\alpha(s) & = \left. \frac{\mathrm d}{\mathrm dz} \tilde{G}^\alpha (z,s)
						\right|_{z=1} \\
						& = s^{\alpha-1} \frac{[-\lambda (z-a_1) (z-a_2) ]\left[(i+1)z^i
						+ a_2^{i+1} (1-a_2)^{-1} \right] + \lambda (2-a_1-a_2)}{\lambda^2(1-a_1)^2(1-a_2)^2} \notag \\
						& = s^{\alpha-1} \left[ - \frac{i+1+a_2^{i+1}(1-a_2)^{-1}}{\lambda (1-a_1)(1-a_2)}
						+ \frac{2-a_1-a_2}{\lambda(1-a_1)^2(1-a_2)^2} \right] \notag \\
						& = s^{\alpha-1} \left[ \frac{i+1+a_2^{i+1}(1-a_2)^{-1}}{s^\alpha}
						+ \frac{2\lambda - (s^\alpha+\lambda+\mu)}{\lambda^2(1-a_1)^2(1-a_2)^2}\right] \notag \\
						& = s^{\alpha-1} \left[ \frac{i+1+a_2^{i+1}(1-a_2^{-1})}{s^\alpha}
						+ \frac{2\lambda - (s^\alpha+\lambda+\mu)}{s^{2\alpha}} \right] \notag \\
						& = s^{\alpha-1} \left[ \frac{1+i +a_2^{i+1}(1-a_2)^{-1}}{s^\alpha}
						+\frac{\lambda-\mu}{s^{2\alpha}} - \frac{s^\alpha}{s^{2\alpha}} \right] \notag \\
						& = s^{\alpha-1} \left[ \frac{i}{s^\alpha} + \frac{\lambda-\mu}{s^{2\alpha}}
						+ \frac{a_2^{i+1}(1-a_2)^{-1}}{s^\alpha} \right] \notag \\
						& = \frac{i}{s} + \frac{\lambda-\mu}{s^{\alpha+1}}
						+ \frac{s^{\alpha-1} a_2^{i+1}(1-a_2)^{-1}}{s^\alpha} \notag \\
						& = \frac{i}{s} + \frac{\lambda-\mu}{s^{\alpha+1}}
						+ \mu\frac{\tilde{p}_0^\alpha(s)}{s^\alpha}. \notag
					\end{align}
					Note the in the above calculation we have used the relation \eqref{rel}.
					Result \eqref{lap-mea}
					immediately yields
					\begin{align}
						\mathbb{E} N^\alpha(t) = i + (\lambda-\mu) \frac{t^\alpha}{\Gamma(\alpha+1)}
						+ \mu J^\alpha p_0^\alpha(t).
					\end{align}
				
				\end{proof}
			
			\end{theorem}			
			
			\begin{remark}			
							
				The validity of Theorem \ref{teo1}
				can be checked with the aid of formula \eqref{meap} as follows.
				\begin{align}
					\mathbb{E} N^\alpha(t) & = \int_0^\infty \mathbb{E} N^1(w) \Pr \{ E^\alpha(t) \in \mathrm dw \} \\
					& = i + (\lambda-\mu) \int_0^\infty w \Pr \{ E^\alpha(t) \in \mathrm dw \}
					+ \mu \int_0^\infty \int_0^t p_0^1(y) \mathrm dy \Pr \{ E^\alpha(t) \in \mathrm dw \}.  \notag
				\end{align}
				Therefore the time-Laplace transform, recalling that $\int_0^\infty e^{-st} \mathrm dt \Pr \{ E^\alpha(t)
				\in \mathrm dw \} = s^{\alpha-1} e^{-ws^\alpha} \mathrm dw$, can be written as
				\begin{align}
					\tilde{\mathbb{E}} N^\alpha(s) & = \frac{i}{s} + \int_0^\infty w s^{\alpha-1} e^{-ws^\alpha}
					\mathrm dw + \mu \int_0^\infty \int_0^w p_0^1(y) \mathrm dy s^{\alpha-1} e^{-ws^\alpha} \mathrm dw \\
					& = \frac{i}{s} + \frac{\lambda-\mu}{s^{\alpha+1}} + \mu s^{\alpha-1}
					\int_0^\infty p_0^1(y) \mathrm dy \int_y^\infty e^{-ws^\alpha} \mathrm dw \notag \\
					& = \frac{i}{s} + \frac{\lambda-\mu}{s^{\alpha+1}} + \mu s^{\alpha-1}
					\int_0^\infty p_0^1(y) \mathrm dy \frac{e^{-y s^\alpha}}{s^\alpha} \notag \\
					& = \frac{i}{s} + \frac{\lambda-\mu}{s^{\alpha+1}} + \mu \frac{\tilde{p}_0^1(s^\alpha)}{s}. \notag
				\end{align}
				Now, by noticing that $\tilde{p}_0^1(s^\alpha)/s = \tilde{p}_0^\alpha(s)/s^\alpha$
				(see formula \eqref{pirelli}) we arrive at
				\begin{align}
					\tilde{\mathbb{E}} N^\alpha(s) = \frac{i}{s} + \frac{\lambda-\mu}{s^{\alpha+1}}
					+ \mu \frac{\tilde{p}_0^\alpha(s)}{s^\alpha},
				\end{align}			
				which leads to formula \eqref{meap}.
				
			\end{remark}

			\begin{remark}
						
				A different form of formula \eqref{lap-mea} can be achieved by writing
				\begin{align}
					\tilde{\mathbb{E}} N^\alpha(s) & = \frac{i}{s} + \frac{\lambda-\mu}{s^{\alpha+1}}
					+ \frac{s^{\alpha-1}a_2^{i+1}}{s^\alpha(1-a_2)} \\
					& = \frac{i}{s} + \frac{\lambda-\mu}{s^{\alpha+1}} - \frac{s^{\alpha-1}
					\lambda a_2^{i+1}(1-a_1)}{s^{2\alpha}},	\notag
				\end{align}
				where we used the fact that $1/(1-a_2) = -\lambda(1-a_1)/s^\alpha$. Furthermore, after considering
				\begin{align}
					a_1 = \frac{s^\alpha + \lambda + \mu}{\lambda} - a_2
					= \frac{\mu}{\lambda a_2},
				\end{align}
				we arrive at
				\begin{align}
					\tilde{\mathbb{E}} N^\alpha(s)
					& = \frac{i}{s} + \frac{\lambda-\mu}{s^{\alpha-1}} + \frac{s^{\alpha-1}a_2^i
					(\mu-\lambda a_2)}{s^{2\alpha}} \\
					& = \frac{i}{s} + \frac{\lambda-\mu}{s^{\alpha-1}} + s^{\alpha-1}
					\frac{\mu a_2^i - \lambda a_2^{i+1}}{s^{2\alpha}} \notag \\
					& = \frac{i}{s} + \frac{\lambda-\mu}{s^{\alpha-1}} +
					\frac{\mu a_2^i - \lambda a_2^{i+1}}{s^{\alpha+1}}. \notag
				\end{align}

			\end{remark}
			
			Let us now address the problem of finding explicit results for the state probabilities
			$p_k^\alpha(t) = \Pr \{ N^\alpha(t) = k | N^\alpha(0) = i \}$ of the proposed fractional queue model.
			We start by using the subordination relation stated in Theorem \eqref{teo1} with
			the classical solution of the M/M/1 queue in terms of modified Bessel functions
			of the first kind. In the non-fractional case ($\alpha=1$) we have \citep[Page 154]{bailey1990elements}
			\begin{align}
				p_k^1(t) = {} & \left( \frac{\lambda}{\mu} \right)^{\frac{1}{2}(k-2)} e^{-(\lambda+\mu)t}
				I_{i-k} \left( 2(\lambda \mu)^{1/2} t \right) \\
				& + \left( \frac{\lambda}{\mu} \right)^{\frac{1}{2}(k-i)}
				\int_0^t e^{-(\lambda +\mu)\tau} \left\{ \lambda I_{i+k+2}
				\left( 2(\lambda \mu)^{1/2} \tau \right) -2(\lambda \mu)^{1/2}
				I_{i+k+1} \left( 2(\lambda \mu)^{1/2} \tau \right) \right. \notag \\
				& \left. + \mu I_{i+k} \left( 2(\lambda \mu)^{1/2} \tau \right) \right\}
				\mathrm d\tau,  \notag
			\end{align}
			where $I_\nu(z)$ is the modified Bessel function of the first kind.
			
			The state probabilities $p_k^\alpha(t)$, $t \ge 0$, $k \ge 0$, $\alpha \in (0,1]$ can thus be determined
			formally by subordination in the following way:
			\begin{align}
				p_k^\alpha(t) = \int_0^\infty p_k^1(y) \Pr \{ E^\alpha(t) \in \mathrm dy \}.
			\end{align}
			Using the time-Laplace transform $\tilde{p}_k^\alpha(s) = \int_0^\infty e^{-st} p_k^\alpha(t)
			\mathrm dt$ we have
			\begin{align}
				\tilde{p}_k^\alpha(s) = {} & \int_0^\infty p_k^1(y) s^{\alpha-1} e^{-ys^\alpha} \mathrm dy \\
				= {} & \left( \frac{\lambda}{\mu} \right)^{\frac{1}{2}(k-2)} \int_0^\infty e^{-(\lambda+\mu)y}
				I_{i-k} \left( 2(\lambda \mu)^{1/2} y \right) s^{\alpha-1} e^{-y s^\alpha} \mathrm dy \notag \\
				& + \left( \frac{\lambda}{\mu} \right)^{\frac{1}{2}(k-i)} \int_0^\infty
				\left[ \int_0^y e^{-(\lambda +\mu)\tau} \left\{ \lambda I_{i+k+2}
				\left( 2(\lambda \mu)^{1/2} \tau \right) -2(\lambda \mu)^{1/2}
				I_{i+k+1} \left( 2(\lambda \mu)^{1/2} \tau \right) \right. \right. \notag \\
				& \left. \left. + \mu I_{i+k} \left( 2(\lambda \mu)^{1/2} \tau \right) \right\}
				\mathrm d\tau \right] s^{\alpha-1} e^{-y s^\alpha} \mathrm dy \notag \\
				= {} & \left( \frac{\lambda}{\mu} \right)^{\frac{1}{2}(k-2)} s^{\alpha-1}
				\int_0^\infty e^{-y(s^\alpha + \lambda + \mu)} I_{i-k} \left( 2(\lambda \mu)^{1/2} y \right) \notag \\
				& + \left( \frac{\lambda}{\mu} \right)^{\frac{1}{2}(k-i)} s^{\alpha-1}
				\int_0^\infty e^{-(\lambda+\mu)\tau} \left\{ \lambda I_{i+k+2}
				\left( 2(\lambda \mu)^{1/2} \tau \right) - 2(\lambda \mu)^{1/2}
				I_{i+k+1} \left( 2(\lambda \mu)^{1/2} \tau \right) \right. \notag \\
				& \left. + \mu I_{i+k} \left( 2(\lambda \mu)^{1/2} \tau \right) \right\} \mathrm d\tau
				\int_\tau^\infty e^{-ys^\alpha} \mathrm dy, \notag
			\end{align}
			Applying the well-known
			Laplace transform for $I_\nu(z)$ we get
			\begin{align}
				\tilde{p}_k^\alpha(s) = {} & \left( \frac{\lambda}{\mu} \right)^{\frac{1}{2}(k-2)} s^{\alpha-1}
				\left[ 2(\lambda \mu)^{1/2} \right]^{k-i} \left[ s^\alpha+\lambda+\mu -\sqrt{(s^\alpha+\lambda+\mu)^2
				-\left[ 2(\lambda \mu)^{1/2} \right]^2} \right]^{i-k} \\
				& \times \left[ (s^\alpha+\lambda+\mu)^2 - \left[ 2(\lambda \mu)^{1/2} \right]^2
				\right]^{-\frac{1}{2}} \notag \\
				& + \left( \frac{\lambda}{\mu} \right)^{\frac{1}{2}(k-i)} \frac{1}{s}
				\left[ \lambda \int_0^\infty e^{-(s^\alpha +\lambda +\mu)\tau} I_{i+k+2}
				\left( 2(\lambda \mu)^{1/2} \tau \right) \mathrm d\tau \right. \notag \\
				& \left. - 2(\lambda\mu)^{1/2} \int_0^\infty e^{-(s^\alpha+\lambda+\mu)\tau}
				I_{i+k+1} \left( 2(\lambda \mu)^{1/2} \tau \right) \right. \notag \\
				& \left. + \mu \int_0^\infty e^{-(s^\alpha +\lambda+\mu)\tau} I_{i+k}
				\left( 2(\lambda \mu)^{1/2} \tau \right) \mathrm d\tau \right] \notag \\
				= {} & \left( \frac{\lambda}{\mu} \right)^{\frac{1}{2}(k-2)}
				s^{\alpha-1} \left[ \frac{s^\alpha + \lambda + \mu
				- \sqrt{(s^\alpha+\lambda+\mu)^2 - 4\lambda\mu}}{2(\lambda\mu)^{1/2}} \right]^{i-k}
				\left[ (s^\alpha+\lambda+\mu)^2 -4\lambda\mu \right]^{-1/2} \notag \\
				& + \left( \frac{\lambda}{\mu} \right)^{\frac{1}{2}(k-i)} \frac{\lambda}{s}
				\left[ 2(\lambda\mu)^{1/2} \right]^{-(i+k+2)}
				\left[ s^\alpha+\lambda+\mu - \sqrt{(s^\alpha+\lambda+\mu)^2-4\lambda\mu} \right]^{i+k+2} \notag \\
				& \times \left[ (s^\alpha+\lambda+\mu)^2 -4\lambda\mu \right]^{-1/2} \notag \\
				& - \left( \frac{\lambda}{\mu} \right)^{\frac{1}{2}(k-i)}
				\frac{2(\lambda\mu)^{1/2}}{s} \left[ 2(\lambda\mu)^{1/2} \right]^{-(i+k+1)}
				\left[ s^\alpha+\lambda+\mu - \sqrt{(s^\alpha+\lambda+\mu)^2 -4\lambda\mu} \right]^{i+k+1} \notag \\
				& \times \left[ (s^\alpha+\lambda+\mu)^2 -4\lambda\mu \right]^{-1/2} \notag \\
				& + \left( \frac{\lambda}{\mu} \right)^{\frac{1}{2}(k-i)} \frac{\mu}{s}
				\left[ 2(\lambda\mu)^{1/2} \right]^{-(i+k)}
				\left[ s^\alpha+\lambda+\mu - \sqrt{(s^\alpha+\lambda+\mu)^2 -4\lambda\mu} \right]^{i+k} \notag \\
				& \times \left[ (s^\alpha+\lambda+\mu)^2 -4\lambda\mu \right]^{-1/2} \notag \\
				= {} & s^{\alpha-1} \frac{(\lambda/\mu)^{\frac{1}{2}(k-2)}}{\sqrt{(s^\alpha+\lambda+\mu)^2
				-4\lambda\mu}} \left[ \frac{s^\alpha+\lambda+\mu-\sqrt{(s^\alpha+\lambda+\mu)^2
				-4\lambda\mu}}{2(\lambda\mu)^{1/2}} \right]^{i-k} \notag \\
				& + s^{\alpha-1} \frac{(\lambda/\mu)^{\frac{1}{2}(k-1)}\lambda}{s^\alpha\sqrt{(s^\alpha+\lambda+\mu)^2
				-4\lambda\mu}} \left[ \frac{s^\alpha+\lambda+\mu-\sqrt{(s^\alpha+\lambda+\mu)^2
				-4\lambda\mu}}{2(\lambda\mu)^{1/2}} \right]^{i+k+2} \notag \\
				& - s^{\alpha-1} \frac{(\lambda/\mu)^{\frac{1}{2}(k-1)}
				2(\lambda\mu)^{1/2}}{s^\alpha\sqrt{(s^\alpha+\lambda+\mu)^2
				-4\lambda\mu}} \left[ \frac{s^\alpha+\lambda+\mu-\sqrt{(s^\alpha+\lambda+\mu)^2
				-4\lambda\mu}}{2(\lambda\mu)^{1/2}} \right]^{i+k+1} \notag \\
				& + s^{\alpha-1} \frac{(\lambda/\mu)^{\frac{1}{2}(k-1)}\mu}{s^\alpha\sqrt{(s^\alpha+\lambda+\mu)^2
				-4\lambda\mu}} \left[ \frac{s^\alpha+\lambda+\mu-\sqrt{(s^\alpha+\lambda+\mu)^2
				-4\lambda\mu}}{2(\lambda\mu)^{1/2}} \right]^{i+k}. \notag
			\end{align}
			Although the obtained Laplace transform $\tilde{p}_k^\alpha(s)$ has a clear structure it cannot
			be inverted in a simple manner. Note anyway that it should be related to the Laplace transform
			of some generalizations of Bessel functions.
			
			In order to obtain more explicit results we must abandon the classical form of
			the state probabilities in terms of Bessel functions. We exploit instead
			a lesser known but certainly more appealing result due to \citet[Chapter 2]{sharma}.
			In particular we refer to equation (2.2.16) at page 17 which we recall here
			for the reader's convenience. Here $\lambda \ne \mu$.
			\begin{align}
				\label{conve}
				p_k^1(t) = {} & \left(1- \frac{\lambda}{\mu}\right) \left( \frac{\lambda}{\mu} \right)^k
				+ e^{-(\lambda+\mu)t} \left( \frac{\lambda}{\mu} \right)^k \sum_{r=0}^\infty 
				\frac{(\lambda t)^r}{r!} \sum_{m=0}^{k+r+i} (r-m) \frac{(\mu t)^{m-1}}{m!} \\
				& + e^{-(\lambda+\mu)t} \sum_{r=0}^\infty (\lambda t)^{k+r-i} (\mu t)^r
				\left( \frac{1}{r!(k+r-i)!} - \frac{1}{(k+r)!(r-i)!} \right). \notag
			\end{align}
			By means of the above formula in the next theorem we derive an explicit
			expression for the state probabilities $p_k^\alpha(t)$, $k \ge 0$, $t \ge 0$, $\alpha \in (0,1]$.
			
			\begin{theorem}
				The state probabilities $p_k^\alpha(t) = \Pr \{ N^\alpha(t) = k | N^\alpha(0) = i \}$,
				$k \ge 0$, $t \ge 0$, $\alpha \in (0,1]$, read
				\begin{align}
					\label{normalmode}
					p_k^\alpha(t) = {} & \left( 1-\frac{\lambda}{\mu} \right) \left( \frac{\lambda}{\mu} \right)^k
					+ \left( \frac{\lambda}{\mu} \right)^k \sum_{r=0}^\infty \sum_{m=0}^{k+r+i}
					\frac{r-m}{r+m} \binom{r+m}{r} \lambda^r \mu^{m-1} t^{\alpha(r+m)-\alpha}
					E_{\alpha,\alpha(r+m)-\alpha+1}^{r+m} \left[ -(\lambda+\mu)t^\alpha \right] \\
					& +\sum_{r=0}^\infty \left[ \binom{k+2r-i}{r} - \binom{k+2r-i}{k+r} \right]
					\lambda^{k+r-i} \mu^r t^{\alpha(k+2r-i)} E_{\alpha,\alpha(k+2r-i)+1}^{k+2r-i+1}
					\left[ -(\lambda+\mu)t^\alpha \right], \notag
				\end{align}
				where
			    \begin{align}
				    E_{\beta,\gamma}^\delta (w) =
				    \sum_{r=0}^\infty \frac{(\delta)_r w^r}{r! \Gamma(\beta r+\gamma)}
				    = \sum_{r=0}^\infty
				    \frac{w^r \Gamma(\delta+r)}{r!\Gamma(\beta r+\gamma)\Gamma(\delta)},
				    \qquad w,\gamma,\beta,\delta \in \mathbb{C}, \: \Re(\beta)>0,
			    \end{align}
			    is the Generalized Mittag--Leffler function \citep{kilbas}.
			    
			    \begin{proof}
					Recurring to Theorem \ref{teo1} we can write for $k \ge 0$, $\lambda \ne \mu$,
					\begin{align}
						p_k^\alpha(t) = {} & \int_0^\infty p_k^1(s) \Pr \{ E^\alpha(t) \in \mathrm ds \} \\
						= {} & \left( 1-\frac{\lambda}{\mu} \right) \left( \frac{\lambda}{\mu} \right)^k
						+ \left( \frac{\lambda}{\mu} \right)^k \sum_{r=0}^\infty
						\sum_{m=0}^{k+r+i} \frac{r-m}{r!m!} \lambda^r \mu^{m-1} \int_0^\infty
						e^{-(\lambda+\mu)s} s^{r+m-1} \Pr \{ E^\alpha(t) \in \mathrm ds \} \notag \\
						& + \sum_{r=0}^\infty \left( \frac{1}{r!(k+r-i)!} - \frac{1}{(k+r)!(r-i)!} \right)
						\lambda^{k+r-i} \mu^r \int_0^\infty e^{-(\lambda+\mu)s} s^{k+2r-i}
						\Pr \{ E^\alpha(t) \in \mathrm ds \}. \notag
					\end{align}
					Applying the Laplace transform to both terms on the right-hand side we obtain
					\begin{align}
						\label{bd}
						\tilde{p}_k^\alpha(z) = {} & \left( 1-\frac{\lambda}{\mu} \right)
						\left( \frac{\lambda}{\mu} \right)^k z^{-1} + \left( \frac{\lambda}{\mu} \right)^k
						\sum_{r=0}^\infty \sum_{m=0}^{k+r+i} \frac{r-m}{r!m!} \lambda^r \mu^{m-1}
						\int_0^\infty e^{-(\lambda+\mu)s} s^{r+m-1} z^{\alpha-1} e^{-sz^\alpha} \mathrm ds \\
						& + \sum_{r=0}^\infty \left( \frac{1}{r!(k+r-i)!} - \frac{1}{(k+r)!(r-i)!} \right)
						\lambda^{k+r-i} \mu^r \int_0^\infty e^{-(\lambda+\mu)s} s^{k+2r-i}
						z^{\alpha-1} e^{-sz^\alpha} \mathrm ds \notag \\
						= {} & \left( 1-\frac{\lambda}{\mu} \right)
						\left( \frac{\lambda}{\mu} \right)^k z^{-1} + \left( \frac{\lambda}{\mu} \right)^k
						\sum_{r=0}^\infty \sum_{m=0}^{k+r+i} \frac{r-m}{r!m!} \lambda^r \mu^{m-1}
						\frac{z^{\alpha-1}}{(z^\alpha+\lambda+\mu)^{r+m}} (r+m-1)! \notag \\
						& + \sum_{r=0}^\infty \left( \frac{1}{r!(k+r-i)!} - \frac{1}{(k+r)!(r-i)!} \right)
						\lambda^{k+r-i} \mu^r \frac{z^{\alpha-1}}{(z^\alpha+\lambda+\mu)^{k+2r-i+1}}
						(k+2r-i)! \notag \\
						= {} & \left( 1-\frac{\lambda}{\mu} \right)
						\left( \frac{\lambda}{\mu} \right)^k z^{-1} + \left( \frac{\lambda}{\mu} \right)^k
						\sum_{r=0}^\infty \sum_{m=0}^{k+r+i} \frac{r-m}{r+m} \binom{r+m}{r}
						\lambda^r \mu^{m-1} \frac{z^{\alpha-1}}{(z^\alpha+\lambda+\mu)^{r+m}} \notag \\
						& + \sum_{r=0}^\infty \left[ \binom{k+2r-i}{r} - \binom{k+2r-i}{k+r} \right]
						\lambda^{k+r-i} \mu^r \frac{z^{\alpha-1}}{(z^\alpha+\lambda+\mu)^{k+2r-i+1}}. \notag
					\end{align}
					To invert equation \eqref{bd} we use the Laplace transform
					(see formula (2.3.24) of \citet{mathai})
					\begin{align}
						\int_0^\infty e^{-zt} t^{\gamma-1} E_{\beta, \gamma}^\delta (wt^\beta) \mathrm dt
						= \frac{z^{\beta \delta -\gamma}}{(z^\beta-w)^\delta},
					\end{align}
					which immediately leads to \eqref{normalmode}.
			    \end{proof}
			    
			\end{theorem}
			
			\begin{remark}
				When $\alpha=1$, formula \eqref{normalmode} becomes the classical solution \eqref{conve}
				because $E_{1,\delta}^\delta (w) = e^w/\Gamma(\delta)$.
			\end{remark}
			
			\begin{remark}
				Result \eqref{normalmode} is particularly interesting because its first addend
				contains the steady-state solution
				\begin{align}
					{}_s p_k^\alpha(t) = \lim_{t \rightarrow \infty} p_k^\alpha(t)
					= \left( 1-\frac{\lambda}{\mu} \right) \left( \frac{\lambda}{\mu} \right)^k, \qquad k \ge 0.
				\end{align}
				Furthermore, it is worth noticing that this geometric
				distribution coincides with that of the classical case
				$\alpha=1$. The whole difference between the fractional and the non-fractional case
				lies in the transient regime.
			\end{remark}

		\section{Path simulation and parameter estimation for the fractional linear birth-death process}

			\label{secsec2}
			We now focus on a related point process which is relatively simpler to treat.
			Let  $ \mathcal{N}(t)$, $t \ge 0 $ be a classical linear birth-death process with $\lambda k >0$, 
			$\mu k >0$  as its birth and death rates, respectively.  Furthermore,
			define $S_k$, $k \in \mathbb{N} \cup \{0\}$,
			as the sojourn time of the process $\mathcal{N}(t)$, $t \ge 0$ in state $k$,
			i.e.,  given that the process is in state $k$,
			$S_k$ is the time until the process leaves that state.
			It is well-known \citep[Section 3.2, Chapter VI]{karlin} that
			\begin{align}
				\Pr \{ S_k \ge t \} = \exp \left[ - (\lambda+\mu) k t \right],
			\end{align}
			and thus,
			\begin{align}
				\Pr \{ S_k \in \mathrm dt \} / \mathrm dt = (\lambda+\mu)k
				\exp \left[ - (\lambda+\mu) k t \right], \qquad t \ge 0.
			\end{align}
			
			Analogous to the preceding section, in order to produce a fractional
			process related to the classical birth-death process it would be natural to
			substitute the unit-order time-derivative in the governing equations
			with a fractional derivative. This has been already carried out
			in \citet{polito}. In the following, we exploit a subordination relation
			similar to that used in Section \ref{secsec} in order to continue the analysis.
			In particular our aim is to develop methods suitable to simulation
			and parameter estimation that will be also applied to the fractional M/M/1 case
			in the last section.
			
			Recall thus that a fractional linear birth-death process  $\mathcal{N}^\alpha(t)$, $t \ge 0$,
			$\alpha \in (0,1]$  satisfies the subordination-relation \citep{polito}
			\begin{align}
				\mathcal{N}^\alpha(t) \overset{\text{d}}{=} N[E^\alpha(t)],
			\end{align}
			where $E^\alpha(t)$ is the right-inverse process to an $\alpha$-stable subordinator
			defined in the previous section.
            Using the above relation, we can easily calculate the distribution
			of the sojourn or holding times $S_k^\alpha$, $k \in \mathbb{N} \cup \{ 0 \}$ for the fractional
			linear birth-death process $\mathcal{N}^\alpha(t)$, $t \ge 0$, as follows.
			\begin{align}
				\label{wtimes}
				\Pr \{ S_k^\alpha \in \mathrm dt \} / \mathrm dt
				& = (\lambda+\mu)k\int_0^\infty \exp\left[-(\lambda+\mu)k s\right]
				\Pr \{ E^\alpha(t) \in \mathrm ds \} \\
				& = (\lambda+\mu)k t^{\alpha-1} E_{\alpha,\alpha} \left[ -(\lambda+\mu)k t^\alpha \right],
				\qquad  t \ge 0,\: k \geq 1, \: 0<\alpha \leq 1. \notag
			\end{align}
			Hence, the holding or sojourn time $S_k^\alpha$ of the fractional linear birth-death process
			$\mathcal{N}^\alpha (t)$  is Mittag--Leffler distributed. 
			
			Another equivalent way to derive the event time distribution above is to
			replace the unit-order derivative  in equation (3.4) of \citet[page 356]{karlin}
			by the Caputo's fractional derivative operator $D^{\alpha}_t $
			used by \cite{polito}. That is,  
			\begin{equation}
				D^{\alpha}_t \Pr \{ S_k^\alpha \ge t \}=  - (\lambda+\mu)k \Pr \{ S_k^\alpha \ge t \},
			\end{equation}
			and solving the above equation, we obtain
			\begin{equation}
				\Pr \{ S_k^\alpha \ge t \}  = E_{\alpha,1} \left[ -(\lambda+\mu)k t^\alpha \right],
			\end{equation}
			which gives equation \eqref{wtimes}.

			An interesting observation is that the birth and death sojourn  times $B_k$ and $D_k$,
			respectively are no longer independent  for $0<\alpha <1$, i.e.,
			\begin{align}
				P( S_k \geq t ) & =   P(\min \lbrace B_k,D_k \rbrace \geq t ) \\
				&=P(B_k>t, D_k>t) \\
				& \neq  P(B_k>t) P( D_k>t),
			\end{align}

			When the process is in state $k$, $k \in \mathbb{N} \cup \{ 0 \}$,  it  transitions
			to the neighboring states $k + 1$ and $k-1$ with probabilities
			$\lambda/(\lambda+\mu)$ and $\mu/(\lambda+\mu)$,  respectively.
			Following  \citet[page 358]{karlin}, a standard procedure to simulate
			trajectories of a fractional linear birth-death process is as follows:

			\vspace{0.1in}

			\noindent ALGORITHM :

			\begin{description}
				\item i) Fix the birth intensity $\lambda$,  the death intensity $\mu$, and the initial
					population size  $\mathcal{N}^\alpha(0) =m$.
				\item ii) Simulate $S_1^\alpha  \stackrel{\text{d}}{=}
					\mathcal{E}_1^{1/ \alpha}T_\alpha$ and $U \stackrel{\text{d}}{=}U(0,1)$.
				\item iii) If $U<\frac{\lambda}{\lambda + \mu}$ then $\mathcal{N}^\alpha(s_1) = m + 1$.
					Otherwise, $\mathcal{N}^\alpha(s_1) = m - 1$. 
				\item iv) Continuing in the same fashion and supposing that the current
					process state is $\mathcal{N}^\alpha(s_{k-1}) = k$,  generate $S_k^\alpha
					\stackrel{\text{d}}{=} \mathcal{E}_k^{1/ \alpha}T_\alpha$ and $U \stackrel{\text{d}}{=}U(0,1)$.
					If $U<\frac{\lambda}{\lambda + \mu}$ then $\mathcal{N}^\alpha(s_k) = k + 1$.
					Otherwise, $\mathcal{N}^\alpha(s_k) = k - 1$.  Repeat iv) until the desired population
					size is achieved or until extinct.
			\end{description}

\begin{figure}[h!t!b!p!]
				\centering
				\includegraphics[height=3in, width=6in]{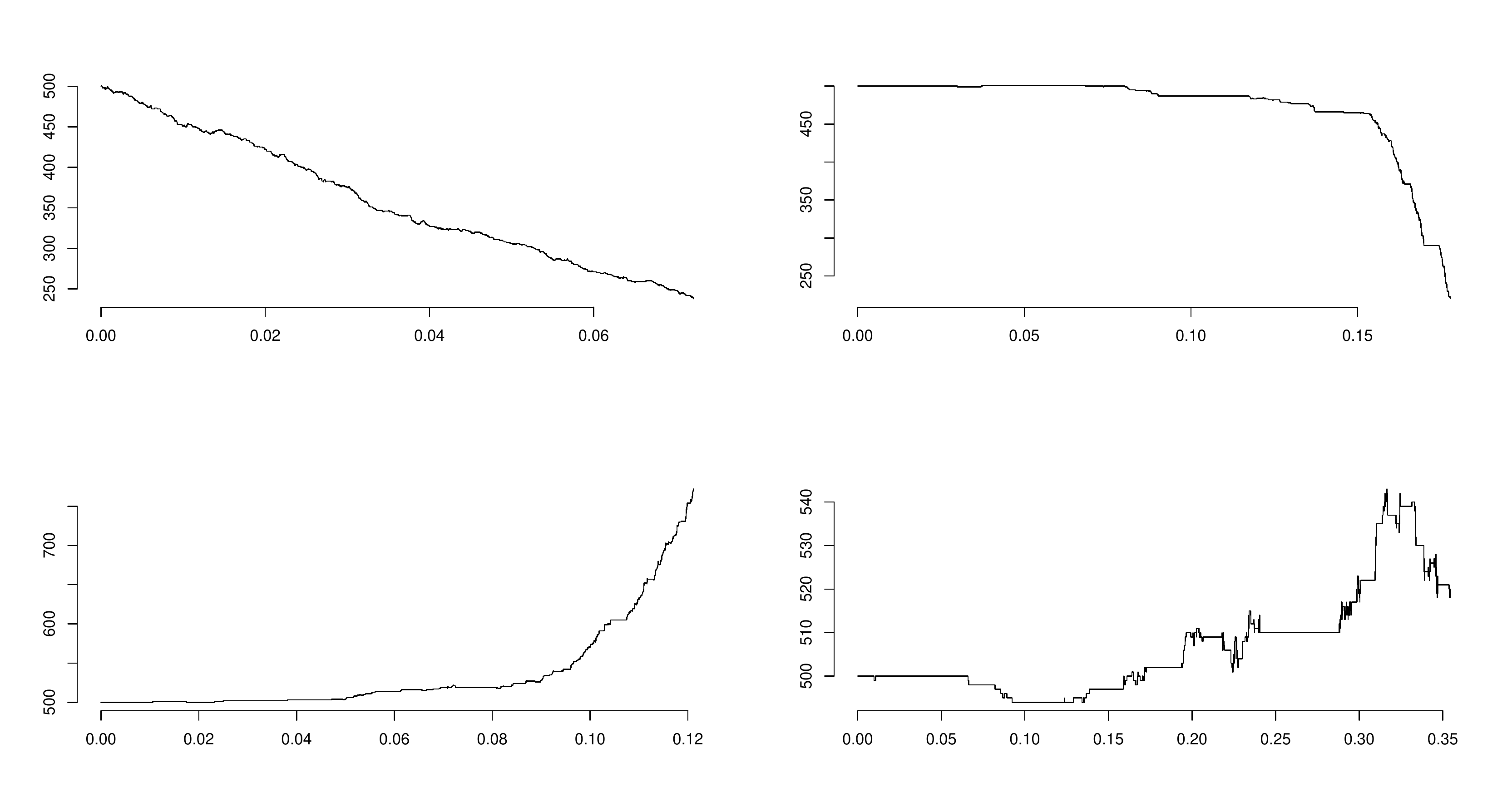}
				\caption{\emph{Sample paths of the fractional  linear birth-death process:  (top left)  $\alpha=1,
					\lambda=5, \mu=15$;  (top right) $\alpha=0.8, \lambda=5, \mu=15$; (bottom left)
					$\alpha=0.8, \lambda=15,
					\mu=5$; (bottom right)  $\alpha=0.8, \lambda= \mu=5$ with an  initial population of $m=500.$}}
				 \renewcommand\belowcaptionskip{0pt}
				\label{Fig1}
			\end{figure}

			Note that $\mathcal{E}_k$ is negative exponentially distributed
			($\mathcal{E}_k \stackrel{\text{d}}{=} \exp \left(  (\lambda + \mu)k \right)$), and
			$T_\alpha$ is a one-sided $\alpha$-stable and independently distributed random variable.
			Below are some sample paths of the	 fractional linear birth-death process.
			The nonlinear trend,  the  longer holding times,  and the slow or  bursting behavior of the
			fractional linear birth-death process   are apparent  in Figure \ref{Fig1}.

			We now provide point estimation algorithms for the parameters $\alpha, \lambda$, and $\mu$.
			Assume that a sample trajectory of size $n$ corresponding to the $n$ random inter-event
			times $S_k^\alpha$'s of the fractional linear birth-death process is observed, where there are $n_B$
			births, $n_D$ deaths, and $n=n_B+n_D$.   Recall the structural representation  of the Mittag--Leffler
			distributed random sojourn time $S_k^\alpha  \stackrel{\text{d}}{=}  \mathcal{E}_k^{1 / \alpha} T_\alpha$,
			where $E_k \stackrel{\text{d}}{=} \exp \left( \theta k \right)$ is independent of a one-sided
			$\alpha^+$-stable distributed random variable  $T_\alpha$, and $\theta=\lambda + \mu$.
			Let $S_k^{\alpha'}=\ln \left( S_k^\alpha \right)$. Then it is well-known that  the mean and variance
			\citep[see details in][]{cuw10} of the log-transformed $k$-th  random  sojourn time
			of the fractional linear birth-death process  are
			\begin{equation}
				\label{assume1}
				\mu_{S_k^{\alpha'}}= \frac{-\ln \left( \theta k \right) }{\alpha} - \gamma,
			\end{equation}
			and
			\begin{equation}
				\label{assume2}
				\sigma_{S_k^{\alpha'}}^2 =\pi^2 \left( \frac{1}{3\alpha^2} -\frac{1}{6} \right),
			\end{equation}
			respectively, where $\gamma \approx 0.5772156649$ is the  Euler--Mascheroni's constant.
			Following \cite{cap12},  the first two moments above therefore suggest that the simple linear
			regression model below can be fitted:
			\begin{equation}
				S_k^{\alpha'} = b_0 + b_1 \ln k + \varepsilon_k, \qquad  k =1,\ldots, n, 
			\end{equation}
			where
			\begin{equation}
				\label{b1e}
				b_0= \frac{-\ln (\theta ) }{\alpha} - \gamma, \qquad b_1= -1/\alpha,
			\end{equation}
			and $\varepsilon_k \stackrel{\text{iid}}{=}  \left(\mu_\varepsilon=0,  \sigma_\varepsilon^2 =
			\sigma_{S_k^{\alpha'}}^2 \right)\stackrel{\text{iid}}{=} \ln \left( \mathcal{E}^{1/ \alpha}
			T_\alpha \right) + \gamma, \mathcal{E} \stackrel{\text{d}}{=} \exp (1)$.
			We point out that  the error distribution depends only on $\alpha$ and is independent of the state $k$
			and $\theta$, and this gives us a simple way of testing the rate fit as follows.
			Generate,  say $m$ samples (each of sample size $n$) from the error distribution using $\widehat{\alpha}$.
			For a fixed significance level, test equality of two parent populations of the observed residuals  and  each
			of the $m$ simulated errors using the two-sample Kolmogorov--Smirnov test, for instance.
			The proportion of the null acceptance out of $m$ tests can then be used to measure model fit. 

			Letting $\sigma_\varepsilon^2 $ or  $\sigma_{S_k^{\alpha'}}^2$ in formula  (\ref{assume2})
			equal to its unbiased estimator
			$\widehat{\sigma}_\varepsilon^2  = \sum_{j=1}^n   \widehat{\varepsilon}_j^2/ (n-2)$,
			we readily obtain the residual-based point estimators
			\begin{equation}
				\label{nuresp}
				\widehat{\alpha}= \left[ {3 \left( \widehat{\sigma}_\varepsilon^2
				\big/ \pi^2 +\frac{1}{6} \right)}  \right]^{-1/2},
				\qquad \qquad
				\widehat{\theta}= \exp\left( -\widehat{\alpha} \left( \widehat{b}_0 + \gamma  \right) \right),
			\end{equation}
			of the model parameters $\alpha$  and $\theta$, correspondingly, where $ \widehat{\varepsilon}_k
			= S_k^{\alpha'} - \widehat{S}_k^{\alpha'}$, and $ \widehat{S}_k^{\alpha'} =\widehat{b}_0 + \widehat{b}_1
			\ln k$. Notice that
			the above estimators exploit the residuals to estimate $\alpha$
			instead of the negative inverse of the least squares (LS) estimate of the slope $b_1$.
			Furthermore, the least squares estimators  of the slope and intercept are
			\begin{equation}
				\widehat{b}_1 = \frac{ \sum_{j=1}^n  S_j^{\alpha'} \left( \ln j - \overline{\ln k }
				\right) }{\sum_{j=1}^n  \left( \ln j -  \overline{\ln k } \right)^2},
				\qquad \qquad
				\widehat{b}_0= \overline{S_k^{\alpha'} }  -  \widehat{b}_1 \cdot \overline{\ln k },
			\end{equation}
			where  $\overline{\ln k } = \sum\limits_{j=1}^n \ln j /n$  and $ \overline{S_k^{\alpha'} }
			=\sum\limits_{j=1}^n S_j^{'} /n$. Hence, the closed-form  point estimators of the intensities
			$\lambda$ and $\mu$ are
			\begin{equation}
				\widehat{\lambda} = \frac{\text{\# of births}}{n} \cdot \widehat{\theta} =\frac{n_B}{n}
				\cdot \widehat{\theta}
			\end{equation}
			and
			\begin{equation}
				\widehat{\mu} = \frac{\text{\# of deaths}}{n} \cdot \widehat{\theta}
				= \frac{n_D}{n} \cdot \widehat{\theta} = \widehat{\theta} - \widehat{\lambda}
				= \frac{n_D}{n} \cdot \widehat{\theta} = \left( 1-\frac{n_B}{n} \right) \cdot \widehat{\theta},
			\end{equation}
			respectively. Table \ref{t1} in the appendix shows some test results based on the percent bias
			\begin{align*}			
				100 \times \frac{| \text{average estimate-parameter value}|}{\text{parameter value}}
			\end{align*}				
			and  the coefficient of variation
			\begin{align*}
				\text{CV} = 100 \times \frac{\text{standard deviation of the estimates}}{\text{average estimate}}
			\end{align*}
			using 1000 simulation runs. Note that we replaced the least squares estimator $\widehat{b}_0$
			by the average of $S_k^{\alpha'} - ( 1 / \widehat{\alpha} ) \ln k$ to improve small sample performance.
			Apparently, the proposed point estimators, especially $\widehat{\alpha}$
			performed relatively well even if the sample size is as small as 100.

			We now provide formulas for the interval estimators of the model parameters. It is worth emphasizing that
			the explicit expressions of the estimators can be utilized to obtain resampling-based interval estimates
			especially for relatively small sample sizes.   It is shown in \cite{cuw10} that 
			\begin{equation}
                                        \label{clta}
				\sqrt{n}\left(\widehat{\alpha}-\alpha\right)  \stackrel{\text{d}}{\longrightarrow}
				N \left(0, \frac{\alpha^2\left( 32-20\alpha^2-\alpha^4\right)}{40} \right).
			\end{equation}
			and a residual-based $(1-\epsilon)100\%$ confidence interval for $\alpha$ directly follows as 
			\begin{equation}
				\label{nuresi}
				\widehat{\alpha} \pm z_{\epsilon/2}
				\sqrt{\frac{\widehat{\alpha}^2\left(32-20\widehat{\alpha}^2-\widehat{\alpha}^4\right)}{40 n} },
			\end{equation}
			where $z_{\epsilon/2}$ is the $(1-\epsilon/2)$th quantile of the standard normal distribution, and
			$0 <\epsilon<1$.  We will now show the asymptotic normality of the estimators
			$\widehat{\lambda}$ and $\widehat{\mu}$.

			\begin{theorem}
				Let $\widehat{p}=n_B/n$ and $p=\lambda/\theta$.  Then 
				\begin{equation}
					\sqrt{n}\left(\widehat{\lambda}- \lambda \right) \stackrel{\text{d}}{\longrightarrow}
					N \left(0, \theta^2 p(1-p) + p^2 \sigma_\theta^2 \right)
				\end{equation}
				as $n \to \infty$ where
				\begin{equation}
					\sigma^2_\theta = e^{ -2\alpha (b_0 + \gamma)}
					\left(b_0 + \gamma \right)^2
					\left[ \frac{\alpha^2\left(32-20\alpha^2 -\alpha^4\right) }{40}
					+ n \alpha^2 \sigma_\varepsilon^2 \left(\frac{1}{n}+\frac{\overline{\ln k}^2}{s} \right) \right] ,
				\end{equation}
				and $s = \sum_{j=1}^n  \left(\ln j -\overline{\ln k} \; \right)^2 $.

				\begin{proof}

					It  can be deduced from \cite{cap12} and the asymptotic property of a Bernoulli/binomial
					sampled proportion that
					\begin{equation}
						\sqrt{n}
						\left(
						\begin{array}{c}
							\widehat{p}-p \\
							\\
							\widehat{\theta} -\theta
						\end{array}
						\right)
						\stackrel{\text{d}}{\longrightarrow}
						N \left(\bm{0}, \bm{\Sigma} \right)
					\end{equation}
					as $n \to \infty$, where the variance-covariance matrix $\bm{\Sigma}$ is defined as
					\begin{equation}
						\bm{\Sigma} = \left(
					    \begin{array}{cc}
				    	    p(1-p) & 0 \\
				    	    & \\
				    	    0 & \sigma^2_\theta\\
					    \end{array}
					    \right).
					\end{equation}

					Invoking a standard result on asymptotic theory,  the two-dimensional Central Limit Theorem
					implies that
					\begin{align}
						\sqrt{n}\big(\textbf{h}(\widehat{\bm{\theta}}_n)-\textbf{h}(\bm{\theta})\big)
						\stackrel{\text{d}}{\longrightarrow} N \left(0,
						\bm{\dot{\textbf{h}}}(\bm{\theta})^{\text{T}}\bm{\Sigma}\bm{\dot{\textbf{h}}}(\bm{\theta})\right),
					\end{align}
					where  $\widehat{\bm{\theta}}_n=(\widehat{p}$,
					$\widehat{\theta})^\text{T},$
					$\bf{h}$ is a mapping from $\mathbb{R}^2 \to\mathbb{R}$,
					$\bm{\dot{\textbf{h}}}(\bf{x})$ is continuous in a neighborhood of $\bm{\theta} \in \mathbb{R}^2$, 
					$\textbf{h}(p,\theta) = p\cdot \theta$, and $\bm{\dot{\textbf{h}}}(p,\theta)
					= \left(\theta,p\right)^{\text{T}}$.
					This concludes the proof.
				\end{proof} 
			\end{theorem}

			\begin{theorem}
				Let $\widehat{q} = 1- \widehat{p}=n_D/n$, $p= \lambda / \theta$, and  $q=\mu / \theta$.
				Then 
				\begin{equation}
					\sqrt{n}\left(\widehat{\mu}- \mu \right) \stackrel{\text{d}}{\longrightarrow}
					N \left(0,\theta^2 p(1-p)^2 + (1- p)^2 \sigma_\theta^2 \right)
				\end{equation}
				as $n \to \infty$.  

				\begin{proof} The proof directly follows from the preceding theorem except that here we consider
					$\textbf{h}(p,\theta) = (1-p)\cdot \theta$ and $\bm{\dot{\textbf{h}}}(p,\theta)=
					\left(- \theta,(1-p) \right)^{\text{T}}$.
				\end{proof} 
			\end{theorem}

			We can now approximate the $(1-\epsilon)100\%$ confidence interval for  $\lambda$ and $\mu$ as
			\begin{equation}
				\widehat{\lambda} \pm z_{\epsilon/2} \cdot \widehat{\sigma}_\lambda, \qquad \qquad  \text{and}
				\qquad  \qquad   \widehat{\mu} \pm z_{\epsilon/2} \cdot \widehat{\sigma}_\mu,
			\end{equation}
			respectively, where
			\begin{equation}
				\widehat{\sigma}_\lambda = \sqrt{\frac{\widehat{\theta}^2 \,\widehat{p} \,  \widehat{q}
				+  \widehat{p}^2 \, \widehat{\sigma}_\theta^2}{n}},
			\end{equation}
			and
			\begin{equation}
				\widehat{\sigma}_\mu = \sqrt{\frac{\widehat{\theta}^2 \, \widehat{p} \, \widehat{q}^2
				+ (1-\widehat{p})^2 \widehat{\sigma}_\theta^2}{n}}.
			\end{equation}

			We now calculate the coverage probabilities using sample sizes $n=10^2, 10^3, 10^4$ and $10^3$
			simulations to test our interval estimators of $\lambda$ and $\mu$ only.
			Notice that the interval estimator for $\nu$ has already been  shown to perform well
			in past related studies \citep[see][]{cuw10,cap12}.
			Table \ref{t2} of the appendix clearly illustrates that the coverage probabilities of the interval estimators
			are closely approaching the true confidence level  when $n$ is at least 1000.
			If  a narrower interval and a larger coverage are preferred then our simulations suggested
			that the previous point estimate replacement  can be used instead. Note  that this replacement  and
			the above simple fit testing schemes can be directly applied to the fractional birth and fractional
			death processes  in \cite{cap12} to enhance performance of the point and interval estimators as well.

			Overall,  Tables \ref{t1}  and \ref{t2}   provide additional merit to the proposed point and interval estimators
			of the model parameters. Aside from the computational simplicity of the proposed parameter estimation
			methods, the  rate fits  can also be checked straightforwardly.

		\section{Trajectory generation and parameter estimation for the fractional
			simple linear birth-death or M/M/1 queue process}

			\label{secsec3}
			Let  $N(t)$, $t \ge 0$ be a classical simple birth-death process with $\lambda>0$, 
			$\mu>0$  as its constant birth and death rates, respectively.
			Furthermore,  define $S_k$, $k \in \mathbb{N} \cup \{0\}$,
			as the sojourn time of the process $N(t)$, $t \ge 0$ in state $k$,
			i.e.\ given that the process is in state $k$,
			$S_k$ is the time until the process leaves that state. Then from the preceding sections,
			it can  easily be deduced  that the holding/sojourn time $S_k^\alpha$'s
			of the fractional simple birth-death or M/M/1 queue
			$N_\alpha (t)$  are independently and identically (IID) Mittag--Leffler distributed, i.e., 
			\begin{align}
				\Pr \{ S_k^\alpha \in \mathrm dt \} / \mathrm dt = (\lambda+\mu) t^{\alpha-1}
				E_{\alpha,\alpha} \left[ - (\lambda+\mu)  t^\alpha \right], \qquad t \ge 0, \: \alpha \in (0,1].
			\end{align}
			Note that everything here immediately follows from the previous section except that
			$\mathcal{E}_k \stackrel{IID}{=} \exp \left(  \theta \right)$.
			Assume that a sample trajectory of population size $n$ corresponding to the $n$ IID random
			inter-event times $S_k^\alpha$'s of the fractional simple birth-death  or M/M/1 process is observed,
			where there are $n_B$ births, $n_D$ deaths, and $n=n_B+n_D$.
			From \cite{cuw10}, a method-of-moments estimator for $\alpha$ is
			\begin{equation}
				\widehat{\alpha}=\frac{\pi}{\sqrt{3\left(\widehat{\sigma}_{S_k^{\alpha'}}^2 +\pi^2/6\right)}}
			\end{equation}
			and
			\begin{equation}
				\widehat{\theta}=\exp\bigg( -\widehat{\alpha}\,\big(
				\widehat{\mathbb{E}}S_k^{\alpha'} + \mathbb{\gamma}\big) \bigg)
			\end{equation}
			is an estimator for $\theta$.  Recall that the asymptotic normality of $\widehat{\alpha}$ follows from the earlier result (\ref{clta}). The appendix's Table \ref{t3} shows some test results based on the percent bias
			and  CV using 1000 simulation runs.  Apparently, the proposed point estimators performed even better than the ones in the linear case.

 			As in the preceding section, we provide formulas for the interval estimators of the model parameters.
			The explicit expressions of the estimators can be used to obtain resampling-based interval estimates
			especially for small sample sizes. A $(1-\epsilon)100\%$ confidence interval for $\alpha$
			is directly obtained from the previous section by simply  replacing the point estimator of $\alpha$. 

			\begin{theorem}
				Let $\widehat{p}=n_B/n$ and $p=\lambda/\theta$.  Then 
				\begin{equation}
					\sqrt{n}\left(\widehat{\lambda}- \lambda \right) \stackrel{\text{d}}{\longrightarrow}
					N \left(0, \theta^2 p(1-p) + p^2 \sigma_\theta^2 \right)
				\end{equation}
				as $n \to \infty$ where
				\begin{equation}
					\sigma^2_\theta = \frac{\theta^2\Big[20\pi^4(2-\alpha^2)-3\pi^2(\alpha^4+20\alpha^2-32)
					(\ln\theta)^2- 720\alpha^3(\ln \theta)\zeta (3)\Big]}{120\pi^2},
				\end{equation}
				and where $\zeta (3)$ is the Riemann-zeta function  evaluated at 3.

				\begin{proof} We omit the routine proof as it follows from the previous theorem.
				\end{proof} 
			\end{theorem}

			\begin{theorem}
				Let $\widehat{q} = 1- \widehat{p}=n_D/n$, $p= \lambda / \theta$, and  $q=\mu / \theta$.  Then 
				\begin{equation}
					\sqrt{n}\left(\widehat{\mu}- \mu \right) \stackrel{\text{d}}{\longrightarrow}
					N \left(0, \theta^2 p(1-p)^2 + (1- p)^2 \sigma_\theta^2 \right)
				\end{equation}
				as $n \to \infty$. 
				\begin{proof} This directly follows from the preceding theorem. 
				\end{proof} 
			\end{theorem}

			We now test our interval estimators for $\lambda$ and $\mu$ by calculating the coverage probabilities using sample sizes
			$n=10^2, 10^3$, and  $10^4$ simulations. 	Table \ref{t4} of the appendix clearly demonstrates that the coverage probabilities of the interval estimators start to 	approach the true confidence level  when $n$ is at least 1000.

%\pagebreak

			In general,  the empirical tests indicate better performance of the proposed point and interval estimators than
			the procedures  for  the fractional  linear  birth-death process.

		\section{Application}

			We demonstrate our methods using two real financial datasets: 1) the monthly Standard \& Poor's (S\&P) index from January 1,  1980 until August 13, 2013 with 248 positive and 155 negative changes; 2) the semi-annual Dow Jones Industrial Average (DJIA) from 1970 until 2013 with 58 positive and  28 negative changes. The data can be downloaded directly from \texttt{finance.yahoo.com} and \texttt{http://www.djindexes.com}, respectively.  In particular, we apply the simple fractional birth-death or M/M/1 queue
			to model the number of positive-negative index changes. 

			Table \ref{t5}  provides the point and  95\% interval estimates for the two financial datasets. 

			\begin{table*}[h!t!b!p!]
				\caption{\emph{Point and 95\% interval estimates.}}
				\vspace{0.1in}
				\centerline{
				\begin{tabular*}{4in}{c||cc|cc}
	& \multicolumn{2}{c|}{S\&P Data}  & \multicolumn{2}{c}{DJIA Data} \\
					Parameter &  Point &  Interval &   Point &  Interval\\
		 \hline   \hline
					& & & & \\
        	        $\alpha$ & 0.949 &  (0.895, 1.002) &  0.897 &  (0.780, 1.014)  \\ 
					& &  & &\\
        	        $ \lambda $  & 0.032 & (0.024,  0.041) & 0.004 & (0.001,  0.008)   \\
					& & & & \\
	         $\mu$     & 0.020  &  (0.015, 0.026)  & 0.002   &  (0.000, 0.004) \\
& & & & \\
\hline
				\end{tabular*}
				}
\label{t5}
			\end{table*}
The above estimates simply suggest that the  monthly S\&P and  semi-annual DJIA changes are highly likely to be non-standard birth-death processes.  We  also examined the rate fit by simulating 1000 samples using the point estimates, and  tested the equality
			of two parent populations  using the two-sample Kolmogorov-Smirnov's test.  The proportions of p-values larger than 0.05 are 98.5\% and 98.6\% for the two datasets, accordingly, which significantly indicate good model fit. 
			These real-world examples clearly demonstrate that the  proposed fractional birth-death model is more general and can also be used as   a	proof-of-concept  or a smoothing tool for the standard  birth-death process.

\section{Appendix}

			\begin{center}
				\begin{table*}[h!t!b!p!]
					\caption{\emph{Percent bias and dispersion of the proposed point estimators
						of  $\alpha, \lambda$ and $\mu$ in Section 3.}} 
					\vspace{0.1in}
					\centerline {
					%\begin{tabular*}{5.5in}{@{\extracolsep{\fill}}|l||l|c@{\hspace{0.01in}}c@{\hspace{0.01in}}|c@{\hspace{0.01in}}c@{\hspace{0.01in}}|c@{\hspace{0.01in}}c@{\hspace{0.01in}}|}
					\begin{tabular*}{5.4in}{@{\extracolsep{\fill}}c|c@{\hspace{0.1in}}||cc|cc|cc}
						\multirow{2}{*}{$(\alpha, \lambda, 
						\mu)$} & \multirow{2}{*}{Estimator}  &  \multicolumn{2}{c|}{$n=10^2$}
						&  \multicolumn{2}{c|}{$n=10^3$} &  \multicolumn{2}{c}{$n=10^4$}  \\
						%\hline Estimator
						&   & Bias & CV& Bias& CV & Bias & CV \\
						\hline \hline
						\multirow{3}{*}{$(0.1, 0.5 , 9)$}  &   $\widehat{\alpha}$     & 0.989 & 8.902 & 0.075 & 2.880   &  0.010  & 0.903 \\ 
			            &  $\widehat{\lambda}$     & 22.389 & 46.798   &2.920  &  25.643 &  0.394 & 10.729   \\
			            &  $\widehat{\mu}$     & 24.938  & 46.696    &3.568 &  23.146  & 0.363 & 9.358  \\
						\hline
						\multirow{3}{*}{$(0.25, 1, 6)$}  &  $\widehat{\alpha}$     & 1.125  & 8.919    &0.019 & 2.667  & 0.028 & 0.909  \\ 
			            &  $\widehat{\lambda}$     & 22.480  & 45.615   &4.749  & 21.624 &  0.020 & 9.309 \\
		 	            &  $\widehat{\mu}$     & 23.273 & 45.997   &  4.309  & 20.478  & 0.171 &  9.467  \\
						\hline
						\multirow{3}{*}{$(0.5, 5, 5)$}  &  $\widehat{\alpha}$     & 0.709  & 8.279    &0.140  & 2.474  & 0.016 & 0.822  \\ 
			            &  $\widehat{\lambda}$     & 20.249  & 42.057    &3.343  & 21.170 &  0.545 &  8.655 \\
			            &  $\widehat{\mu}$     & 20.992 & 41.142    &  3.767  & 21.090  & 0.621 & 8.637 \\
						\hline
						\multirow{3}{*}{$(0.75, 7, 1)$}  & $\widehat{\alpha}$     & 0.316  & 7.213    &0.109  & 2.272  & 0.004 & 0.679  \\ 
		 	            &  $\widehat{\lambda}$     & 10.372  & 38.858   &1.945  & 17.192 &  0.286 & 6.991 \\
				        &  $\widehat{\mu}$     & 10.205  & 42.200   &2.443 & 19.205  & 0.204 & 7.472  \\
						\hline
						\multirow{3}{*}{$(0.95, 10, 0.5)$}  &  $\widehat{\alpha}$     & 0.924  & 5.386    &0.077  & 1.764 & 0.008 & 0.532  \\ 
			            &  $\widehat{\lambda}$     & 9.544  & 28.672    &1.539  & 13.572 &  0.272 & 5.186  \\
            			&  $\widehat{\mu}$     &  8.875  & 41.985    &1.559  & 19.671  & 0.468 & 7.193  \\
						\hline
					\end{tabular*}
					}
					\label{t1}
				\end{table*}
			\end{center}

			\begin{center}
				\begin{table*}[h!t!b!p!]
					\caption{\emph{Coverage probabilities of the proposed interval estimators of  $\lambda$ and $\mu$
						using a $95\%$  confidence level  in Section 3.}} 
					\vspace{0.1in}
					\centerline {
					%\begin{tabular*}{5.5in}{@{\extracolsep{\fill}}|l||l|c@{\hspace{0.01in}}c@{\hspace{0.01in}}|c@{\hspace{0.01in}}c@{\hspace{0.01in}}|c@{\hspace{0.01in}}c@{\hspace{0.01in}}|}
					\begin{tabular*}{4in}{@{\extracolsep{\fill}}c|c@{\hspace{0.1in}}||c|c|c}
						$(\alpha, \lambda, \mu)$ & Parameter &  $n=10^2$ & $n=10^3$ & $n=10^4$  \\
						\hline \hline
						\multirow{2}{*}{$(0.1, 10 , 90)$}   &  $\lambda$     & 0.879 & 0.936   &0.954\\
        		        &  $\mu$     & 0.888 & 0.937    & 0.955  \\
						\hline
						\multirow{2}{*}{$(0.25, 70, 30)$}   &  $\lambda$     & 0.892 & 0.925   & 0.955  \\
        		        &  $\mu$     & 0.876 &  0.926  & 0.958 \\
						\hline
						\multirow{2}{*}{$(0.5, 50, 50)$}     &  $\lambda$     & 0.896 & 0.934   & 0.946  \\
        	    	    &  $\mu$     & 0.894  & 0.933   & 0.947  \\
						\hline
						\multirow{2}{*}{$(0.7, 5, 95)$}     &  $\lambda$     & 0.864 & 0.922   & 0.947\\
            		    &  $\mu$     & 0.882 &  0.924   &0.950 \\
						\hline
						\multirow{2}{*}{$(0.95, 20, 80)$}    &  $\lambda$     & 0.900 & 0.948   & 0.951  \\
            		    &  $\mu$     & 0.910  & 0.950   & 0.959 \\
						\hline
					\end{tabular*}
					}
					\label{t2}
				\end{table*}
			\end{center}

			\begin{center}
				\begin{table*}[h!t!b!p!]
					\caption{\emph{Percent bias and dispersion of the proposed point estimators of
					$\alpha, \lambda$ and $\mu$  in Section 4.}} 
					\vspace{0.1in}
					\centerline {
					%\begin{tabular*}{5.5in}{@{\extracolsep{\fill}}|l||l|c@{\hspace{0.01in}}c@{\hspace{0.01in}}|c@{\hspace{0.01in}}c@{\hspace{0.01in}}|c@{\hspace{0.01in}}c@{\hspace{0.01in}}|}
					\begin{tabular*}{5.4in}{@{\extracolsep{\fill}}c|c@{\hspace{0.1in}}||cc|cc|cc}
						\multirow{2}{*}{$(\alpha, \lambda, 
						\mu)$} & \multirow{2}{*}{Estimator}  &  \multicolumn{2}{c|}{$n=10^2$}  &  \multicolumn{2}{c|}{$n=10^3$} &  \multicolumn{2}{c}{$n=10^4$}  \\
						%\hline Estimator
						&   & Bias & CV& Bias& CV & Bias & CV \\
						\hline \hline
						\multirow{3}{*}{$(0.1, 0.5 , 9)$}  &   $\widehat{\alpha}$     & 1.180 & 9.077 & 0.262 & 2.924   &  0.015  & 0.864 \\ 
		       	    	&  $\widehat{\lambda}$     & 4.551 & 45.894   &0.731 &  15.929 &  0.072 & 4.804   \\
        			    &  $\widehat{\mu}$     & 6.017  & 25.210    & 0.964 &  8.649  & 0.122 & 2.728  \\
						\hline
						\multirow{3}{*}{$(0.25, 1, 6)$}  &  $\widehat{\alpha}$     & 0.973  & 8.836    &0.044 & 2.786  & 0.022 & 0.916  \\ 
				        &  $\widehat{\lambda}$     &7.344  & 31.717   &0.270  & 11.330  &  0.054 & 3.501 \\
			            &  $\widehat{\mu}$     & 5.322 & 22.923  &  0.337  & 7.979 & 0.021 &  2.450  \\
						\hline
						\multirow{3}{*}{$(0.5, 5, 5)$}  &  $\widehat{\alpha}$     & 0.506  & 8.043   &0.016  & 2.534  & 0.044 & 0.818  \\ 
			            &  $\widehat{\lambda}$     & 6.411  & 24.339    &0.088  & 8.609 &  0.050 &  2.640 \\
	        		    &  $\widehat{\mu}$     &  5.242 & 28.822    &  0.154 & 8.869 & 0.048 &  2.682 \\
						\hline
						\multirow{3}{*}{$(0.75, 7, 1)$}  & $\widehat{\alpha}$     & 0.793  & 7.568    &0.054  & 2.388  & 0.027 & 0.701  \\ 
			            &  $\widehat{\lambda}$     & 3.933  & 21.190  & 0.452  & 6.796 &  0.204 & 1.901 \\
	        		    &  $\widehat{\mu}$     & 6.330  & 32.114   & 0.497 & 10.516  & 0.098 & 3.254  \\
						\hline
						\multirow{3}{*}{$(0.95, 10, 0.5)$}  &  $\widehat{\alpha}$     & 0.541  & 5.587    &0.019  & 1.795 & 0.007 & 0.558  \\ 
			            &  $\widehat{\lambda}$     & 1.820  & 13.947    &0.030  & 4.508 &  0.031 & 1.540  \\
			            &  $\widehat{\mu}$     &  0.061  & 43.363   &0.147  & 14.001  & 0.121 & 4.539  \\
						\hline
					\end{tabular*}
					}
					\label{t3}
				\end{table*}
			\end{center}
			\begin{center}
				\begin{table*}[h!t!b!p!]
					\caption{\emph{Coverage probabilities of the proposed interval estimators of  $\lambda$ and $\mu$
						using a $95\%$  confidence level  in Section 4.}} 
					\vspace{0.1in}
					\centerline {
%\begin{tabular*}{5.5in}{@{\extracolsep{\fill}}|l||l|c@{\hspace{0.01in}}c@{\hspace{0.01in}}|c@{\hspace{0.01in}}c@{\hspace{0.01in}}|c@{\hspace{0.01in}}c@{\hspace{0.01in}}|}
					\begin{tabular*}{4in}{@{\extracolsep{\fill}}c|c@{\hspace{0.1in}}||c|c|c}
						$(\alpha, \lambda, \mu)$ & Parameter &  $n=10^2$ & $n=10^3$ & $n=10^4$  \\
						\hline \hline
						\multirow{2}{*}{$(0.1, 10 , 90)$}   &  $\lambda$     & 0.921 & 0.950   &0.951\\
                        &  $\mu$     & 0.932 & 0.959   & 0.948 \\
						\hline
						\multirow{2}{*}{$(0.25, 70, 30)$}   &  $\lambda$     & 0.936 & 0.954   & 0.955  \\
                        &  $\mu$     & 0.927 &  0.942  & 0.958 \\
						\hline
						\multirow{2}{*}{$(0.5, 50, 50)$}     &  $\lambda$     & 0.932 & 0.957   & 0.949  \\
                        &  $\mu$     & 0.925  & 0.948   & 0.953 \\
						\hline
						\multirow{2}{*}{$(0.7, 5, 95)$}     &  $\lambda$     & 0.869 & 0.944   & 0.947\\
                        &  $\mu$     & 0.931 &  0.955   &0.950 \\
						\hline
						\multirow{2}{*}{$(0.95, 20, 80)$}    &  $\lambda$     & 0.927 & 0.961   & 0.959  \\
                        &  $\mu$     & 0.917  & 0.947   & 0.947 \\
						\hline
					\end{tabular*}
					}
					\label{t4}
				\end{table*}
			\end{center}

		\section*{Aknowledgement}
		
			The authors are grateful to the reviewers and the editors for significantly improving the paper.  Federico Polito has been supported by project AMALFI (Universit\`{a} di Torino/Compagnia di San Paolo). Dexter Cahoy is also supported  by  Louisiana Board of Regents Research Competitiveness Subprogram grant LEQSF(2011-14)-RD-A-15.

\end{document}